\title{On Idempotent Generated Semigroups}
\author{Jo\~ao Ara\'ujo}
\date{}

\documentstyle[12pt]{article}
\begin{document}
\newtheorem{theorem}{\indent \bf Theorem}
\newtheorem{lemma}[theorem]{\indent \bf Lemma}
\newtheorem{proposition}[theorem]{\indent \bf Proposition}
\newtheorem{definition}{\indent \bf Definition}
\newtheorem{corollary}[theorem]{\indent \bf Corollary}
\newtheorem{remark}[theorem]{\indent \bf Remark}
\newtheorem{problem}[theorem]{\indent \bf Problem}
\newtheorem{conjecture}[theorem]{\indent \bf Conjecture}
\def\diag{\mathop{\rm diag}}
\def\rank{\mathop{\rm rank}}
\def\Ker{\mathop{\rm Ker}}
\def\End{\mathop{\rm End}}
\def\Aut{\mathop{\rm Aut}}
\def\GL{\mathop{\rm GL}}
\def\S{\mathop{\cal S}}
\def\halmos{\hfill\rule{6pt}{6pt}}
\def\R{\mathord{\sf{l\hspace{-0.1em}R}}}
\def\C{\mathord{\hspace{0.15em}\sf{l\hspace{-0.4em}C}}}
\def\Q{\mathord{\hspace{0.15em}\sf{l\hspace{-0.4em}Q}}}
\def\Z{\mathord{\sf{Z\hspace{-0.4em}Z}}}
\def\N{\mathord{\sf{l\hspace{-0.1em}N}}}
\def\D{\mathord{\sf{l\hspace{-0.1em}D}}}
\def\P{\mathord{\sf{l\hspace{-0.1em}P}}}
\def\F{\mathord{\sf{l\hspace{-0.1em}F}}}

 \newcommand{\la}{\langle}
  \newcommand{\ra}{\rangle}

\maketitle

\begin{abstract}
We provide short and direct proofs for some classical theorems proved by  Howie, Levi and McFadden    concerning
idempotent generated semigroups of transformations on a finite set.\

\end{abstract}

Let $n$ be a natural number and $[n]=\{1, \ldots ,n\}$. Let $T_n, S_n$ be, respectively, the transformation
semigroup  and the symmetric group on $[n]$. Let $P=(A_i)_{i\in [k]}$ be a partition of $[n]$ and let
$C=\{x_1,\ldots , x_k\}$ be a cross-section of $P$ (say $x_i\in A_i$). Then we represent $A_i$ by $[x_i]_{_P}$ and
the pair $(P,C)$ induces an idempotent mapping defined by $[x_i]_{_P}e=\{x_i\}$. Conversely, every idempotent can
be so constructed. To save space instead of  $e=\left(\begin{array}{cccccc}
[x_1]_{_P}&\ldots &[x_k]_{_P}\\
x_1      &\ldots & x_k
\end{array}\right)$ we write $e=([x_1]_{_P}, \ldots,[x_k]_{_P})$.
This notation extends to $e=([\underline{x_1},y]_{_P},[x_2]_{_P}, \ldots,[x_k]_{_P})$ when  $y\in [x_1]_{_P}$ and
$[x_i]_{_P}e=\{x_i\}$. By $([x_1],\ldots ,[\underline{x_i},{y}], \ldots,[x_k])$ we denote the set of all
idempotents $e$ with image $\{x_1,\ldots , x_n\}$ and such that the $Ker(e)$-class of $x_i$ contains (at least)
two elements: $x_i$ and $y$, where the underlined element (in this case $x_i$) is the image of the class under
$e$. For $a\in T_n$, we denote the image of $a$ by $\nabla a$.
\begin{lemma}\label{1}Let  $a\in I_n=T_n\setminus S_n$, $rank(a)=k$,  and
 $(xy)\in Sym(X)$. Then   $a(xy)=ab$, where $b=1$ or $b=e_1e_2e_3$, with $e_i^2=e_i$ and $rank(e_i)=k$.
\end{lemma}
\it Proof. \rm Let   $\nabla a=\{a_i \mid i\in [k]\}$. If $x,y\not\in \nabla a$, then $a(xy)=a$.
 If, say, $x=a_1$ and $y\not\in \nabla a$, then $ae_1=a(a_1y)$, for all $e_1\in
 ([a_1,\underline{y}], [a_2], \ldots, [a_k])$. Finally, if $x,y\in \nabla a$, without loss of generality, we can assume  that
   $x=a_1$ and $y=a_{2}$. Then $a(a_1\ a_{2})=ae_2e_3e_4$, for $e_2\in ([a_1],
[a_{2},\underline{u}],[a_3], \ldots ,[a_k])$, $e_3\in ([u], [a_{1},\underline{a_{2}}],[a_3], \ldots ,[a_k])$ and
$e_4 \in ([u,\underline{a_1}], [a_{2}], \ldots ,[a_k])$. \halmos
\begin{theorem}  \cite{howie}
Every ideal of $I_n$ is generated by its own idempotents.
\end{theorem}
\it Proof. \rm Let $a\in I_n$. Then $a=eg$ for some $e=e^2\in T_n$ and $g\in S_n$. Therefore $a=e(x_1y_1)\ldots
(x_m y_m)$ and hence, applying $m$ times Lemma \ref{1}, $a$ can be obtained as a product of idempotents of the
same rank as $a$.\halmos

Let $t\in I_n$ and $g\in S_n$. Denote $g^{-1}tg$ by $t^g$ and let $C_t=\{t^g\mid g\in S_n\}$ and   $t^{S_n}=\la
\{t\}\cup S_n\ra\setminus S_n$. For  $f=([x_1,\underline{w}]_Q,[x_2]_Q,\ldots , [x_k]_Q)$ and $g\in S_n$, it is
easy to check that we have $ f^g=([x_1g,\underline{wg}]_{Qg},\ldots , [x_kg]_{Qg})\in
([x_1g,\underline{wg}],\ldots , [x_kg])$.

% We can sharpen Lemma \ref{1} as follows.
\begin{lemma}\label{11}Let  $a\in I_n$, $rank (a)= k$, let
$f=([x_1,\underline{w}]_Q, [x_2]_Q, \ldots , [x_k]_Q)$ and let
 $(xy)\in S_n$. Then   $a(xy)=ab$, where $b=1$ or $b=e_1e_2e_3$, with $e_i \in C_f$.
\end{lemma}
\it Proof. \rm  As in  Lemma \ref{1}, one only has to show that the sets  $([a_1,\underline{y}], [a_2], \ldots,
[a_k])$, $([a_1], [a_{2},\underline{u}], \ldots ,[a_k])$, $([u], [a_{1},\underline{a_{2}}], \ldots ,[a_k])$ and
$([\underline{a_1},u], [a_{2}], \ldots ,[a_k])$ intersect $C_f$.  Let $g\in S_n$ such that $x_1g=a_1, wg=y$ and
$x_ig=a_i$ ($2\leq i\leq k$). Thus $ f^g=([x_1g,\underline{w}g]_{Qg}, [x_2g]_{Qg}, \ldots , [x_kg]_{Qg})\in
([a_1,\underline{y}], [a_2], \ldots, [a_k])$. The proof that  $C_f$ intersects the remaining three sets is
similar.  \halmos

Repeating the arguments of Theorem \ref{1} together with Lemma \ref{11} we have

\begin{corollary}\label{tres}   Let $eg\in I_n$ (for some $e=e^2$, $g\in S_n$) and let $f^2=f$ with $rank(eg)=rank(f)$.
 Then $eg\in e\la {C_f}\ra$ and, in particular,  $eg\in \la {C_e}\ra$.
\end{corollary}
From now on let $a=eg\in I_n$ (for some $e=e^2$, $g\in S_n$).
\begin{corollary} \label{quatro} $\la {C_e}\ra =e^{S_n}={(eg)}^{S_n}=a^{S_n}$.
\end{corollary}
\it Proof. \rm The only non-trivial inclusion is $e^{S_n} \subseteq \la {C_e}\ra$. As $e^{S_n}$ is generated by
$\{geh\mid g,h \in S_n\}$, let $g,h\in S_n$. By  Corollary \ref{tres}, $he=(heh^{-1}) h \in \la {C_{heh^{-1}}}\ra
$ and $eg\in \la {C_e}\ra $. Since $heh^{-1}\in \la {C_e}\ra $, it follows that $\la {C_{heh^{-1}}}\ra \leq \la
{C_e}\ra$. Thus $heg=(he)(eg) \in \la {C_e}\ra $. It is proved that $e^{S_n} \subseteq \la {C_e}\ra $. \halmos
\begin{theorem}\cite{levi}
  $\la {C_a}\ra =\la {C_e}\ra =a^{S_n}$ and hence is idempotent generated.
\end{theorem}
\it Proof. \rm Since $a=eg$, then $e=ag^{-1}$ and it is easy to check that we have $(ag^{-1})^n=aa^ga^{g^2}\ldots
a^{g^{n-1}}g^{-n}$. But for some natural $m$ we have $g^{-m}=1$. Thus $(ag^{-1})^m=aa^ga^{g^2}\ldots
a^{g^{m-1}}\in \la {C_a}\ra$ so that $e=e^m=(ag^{-1})^m\in \la {C_a}\ra$ and hence, by Corollary \ref{quatro}, it
follows that $a^{S_n}=\la {C_e}\ra \leq \la {C_a}\ra \leq a^{S_n}$. (Cf. \cite{mcalister}). \halmos


\footnotesize

\begin{thebibliography}{99}

\bibitem{howie}
J.\ M.\ Howie, R.\ B.\ McFadden
{\em Idempotent rank in finite full transformation semigroups},
Proc. R. Soc. Edinburgh, Sect.\ A   {\bf 114} (1990), 161--167.

\bibitem{levi}
  I. Levi and R. McFadden  {\em $S_n$-normal semigroups}, Proc. Edinburgh Math. Soc., {\bf 37} (1994), 471-476.


\bibitem{mcalister}
  D.B. McAlister  {\em Semigroups generated by a group and an idempotent},
 Comm. Algebra, {\bf 26} (1998), 515-547.


\end{thebibliography}
\end{document}